\documentclass[11pt]{article}
\usepackage{amsmath}
\usepackage{amsfonts}
\usepackage{amssymb}
\usepackage{dsfont}
\usepackage{hyperref}
\usepackage{indentfirst}
\usepackage{mathrsfs}

\usepackage[top=1in,bottom=1in,left=1.25in,right=1.25in]{geometry}
\date{}
\allowdisplaybreaks
\begin{document}

\renewcommand{\baselinestretch}{1.2}
\renewcommand{\arraystretch}{1.0}

\title{\bf Integral Dual of some infinite dimensional Hopf quasigroups}
\author
{
 \textbf{Tao Yang} \footnote{Corresponding author.
 College of Science, Nanjing Agricultural University, Nanjing 210095, China. E-mail: tao.yang@njau.edu.cn}
}
\maketitle

\begin{center}
\begin{minipage}{12.cm}

 \textbf{Abstract}
 For an infinite dimensional Hopf quasigroup, if the faithful integral exists,
 then it is unique up to scalar. Base on the faithful integrals, we  construct the integral dual of a class of infinite dimensional Hopf quasigroups, and show that the integral dual has a similar structure to Hopf coquasigroup, which is a regular multiplier Hopf coquasigroup with a faithful integral.
\\

\\


\end{minipage}
\end{center}
\normalsize

\section{Introduction}
\def\theequation{\thesection.\arabic{equation}}
\setcounter{equation}{0}

Just as many Lie groups have an entirely algebraic description as commutative Hopf algebras, J. Kim and S. Majid developed the corresponding theory of 'algebraic quasigroups' including the coordinate algebra $k[S^7]$ of the 7-sphere  in \cite{KM}. They defined the notion of a Hopf quasigroup and showed that a theory similar to that of Hopf algebras was possible in this case.
The first author in his following paper \cite{K} developed the integral theory for Hopf (co)quasigroups, Fourier transform, and showed that a finite dimensional Hopf (co)quasigroup has a unique integration up to scalar and an invertible antipode.
Then a natural question arise:
(Q1) Is the integral also unique on an infinite dimensional Hopf quasigroup?

 As shown in \cite{KM},  the dual of a  finite dimensional Hopf quasigroup is a Hopf coquasigroup. Then
 (Q2) How about the dual of an infinite dimensional Hopf quasigroup?
 These two questions motivated this paper.

 In \cite{V98}, A. Van Daele developed the theory of algebraic quantum groups. As an application, it provided  the 'integral dual' of an infinite dimensional Hopf algebra with integrals.
 The non-degenerate faithful integrals play a key role in the dual.
 Inspired by this, we first consider the properties of integrals on the infinite dimensional Hopf quasigroup, and then construct the integral dual. we find that the integral dual has a similar structure of Hopf coquasigoup, which is a multiplier  Hopf coquasigoup of discrete type.

This paper is organized as follows. In section 2, we introduce some notions: Hopf (co)quasigroups, which will be used in the following sections.

In section 3, we consider the integral on an infinite dimensional Hopf quasigroup, and show that the faithful integrals are unique up to scalar.

In section 4, for an infinite dimensional Hopf quasigroup $H$ with a faithful integral $\varphi$ , we construct the integral dual $\widehat{H} = \varphi(\cdot H)$, and show that $\widehat{H}$ is a regular multiplier Hopf coquasigroup with a faithful integral.

In section 5, we introduce our motivating example, and make some comments on multiplier Hopf coquasigroups.

\section{Preliminaries}
\def\theequation{\thesection.\arabic{equation}}
\setcounter{equation}{0}

 Throughout this paper, all spaces we considered are over a fixed field $k$ (e.g., the complex number field $\mathds{C}$).

\subsection{Hopf (co)quasigroup}

 Recall from \cite{KM} a \emph{Hopf quasigroup} is possibly non-associative but unital algbera $H$ equipped with algebra homomorphisms $\Delta: H\longrightarrow H\otimes H$,
 $\varepsilon:H\longrightarrow k$ forming a coassociative coalgebra and a map $S: H\longrightarrow H$ such that
 \begin{eqnarray*}
 && m(id\otimes m)(S\otimes id \otimes id)(\Delta \otimes id) = \varepsilon\otimes id = m(id\otimes m)(id\otimes S \otimes id)(\Delta \otimes id),\\
 && m(m\otimes id)(id\otimes S \otimes id)(id \otimes \Delta) = id\otimes \varepsilon =  m(m\otimes id)(id\otimes id \otimes S)(id \otimes \Delta).
 \end{eqnarray*}
 These two equations can be written more explicitly as: for all $g, h\in H$,
 \begin{eqnarray*}
 \sum S(h_{(1)})(h_{(2)}g) = \sum h_{(1)}\big(S(h_{(2)})g \big) = \sum \big(gS(h_{(1)}) \big)h_{(2)} = \sum (gh_{(1)})S(h_{(2)}) = \varepsilon(h)g,
 \end{eqnarray*}
 where we write $\Delta h = \sum h_{(1)} \otimes h_{(2)}$ and for brevity, we shall omit the summation signs.

 The Hopf quasigroup $H$ is called \emph{flexible} if
 \begin{eqnarray}
 h_{(1)} (gh_{(2)}) = (h_{(1)}g) h_{(2)}, \qquad \forall g, h\in H, \label{2.1}
 \end{eqnarray}
 and \emph{alternative} if also
 \begin{eqnarray}
 h_{(1)} (h_{(2)}g) = (h_{(1)}h_{(2)}) g, \qquad h(g_{(1)} g_{(2)}) = (hg_{(1)}) g_{(2)}, \qquad \forall g, h\in H. \label{2.2}
 \end{eqnarray}
 $H$ is called \emph{Moufang} if
 \begin{eqnarray}
 h_{(1)} \big(g (h_{(2)}f) \big) = \big( (h_{(1)}g) h_{(2)} \big)f, \qquad \forall h, g, f\in H. \label{2.3}
 \end{eqnarray}

 It was proved that the antipode $S$ is antimultiplicative and anticomultiplicative, i.e., for all $g, h\in H$,
 \begin{eqnarray*}
 S(gh) =S(h)S(g), \qquad \Delta(Sh) = S(h_{(2)}) \otimes S(h_{(1)}).
 \end{eqnarray*}
 Moreover, if $H$ is cocommutative flexible Hopf quasigroup, then $S^{2}=id$ and for all $g, h\in H$,
 \begin{eqnarray*}
 h_{(1)}\big( g S(h_{(2)}) \big) = (h_{(1)}g) S(h_{(2)}).
 \end{eqnarray*}

 Dually, we can obtain a Hopf coquasigroup by reversing the arrows on each map in Hopf quasigroup.

 A \emph{Hopf coquasigroup} is a unital associative algebra $A$ equipped with counital algebra homomorphism $\Delta: A\longrightarrow A\otimes A$,
 $\varepsilon: A\longrightarrow k$ and a linear map $S: A \longrightarrow A$ such that for all $a\in A$,
 \begin{eqnarray*}
 && (m\otimes id)(S\otimes id \otimes id)(id \otimes \Delta)\Delta = 1 \otimes id = (m\otimes id)(id\otimes S \otimes id)(id \otimes\Delta)\Delta,\\
 && (id\otimes m)(id\otimes S \otimes id)(\Delta \otimes id)\Delta = id \otimes 1 = (id\otimes m)(id\otimes id \otimes S)(\Delta \otimes id)\Delta.
 \end{eqnarray*}
 In other word,
 \begin{eqnarray*}
 &&S(a_{(1)})a_{(2)(1)} \otimes a_{(2)(2)} = 1 \otimes a = a_{(1)}S(a_{(2)(1)}) \otimes a_{(2)(2)},\\
 && a_{(1)(1)} \otimes S(a_{(1)(2)})a_{(2)} = a \otimes 1 = a_{(1)(1)} \otimes a_{(1)(2)} S(a_{(2)}).
 \end{eqnarray*}

 A Hopf coquasigroup is \emph{flexible} if
 \begin{eqnarray}
 a_{(1)}a_{(2)(2)}\otimes a_{(2)(1)} = a_{(1)(1)}a_{(2)} \otimes a_{(1)(2)}, \qquad \forall a\in A,
 \end{eqnarray}
 and \emph{alternative} if also
 \begin{eqnarray}
 && a_{(1)}a_{(2)(1)}\otimes a_{(2)(2)} = a_{(1)(1)}a_{(1)(2)} \otimes a_{(2)}, \\
 && a_{(1)}\otimes a_{(2)(1)}a_{(2)(2)} = a_{(1)(1)} \otimes a_{(1)(2)}a_{(2)}, \qquad \forall a\in A.
 \end{eqnarray}
 $A$ is called \emph{Moufang} if
 \begin{eqnarray}
 a_{(1)}a_{(2)(2)(1)}\otimes a_{(2)(1)} \otimes a_{(2)(2)(2)} = a_{(1)(1)(1)}a_{(1)(2)} \otimes a_{(1)(1)(2)}\otimes a_{(2)}, \qquad \forall a\in A,
 \end{eqnarray}
 The term 'counital' here means
 \begin{eqnarray*}
 (id\otimes \varepsilon)\Delta = id = (\varepsilon\otimes id)\Delta.
 \end{eqnarray*}
 However, $\Delta$ is not assumed to be coassociative.

 It was shown in \cite{KM} Proposition 5.2 that: Let $A$ be a Hopf coquasigroup, then
 \begin{enumerate}
   \item[(1)] $m(S\otimes id)\Delta = \mu\varepsilon = m(id\otimes S)\Delta$.
   \item[(2)] $S$ is antimultiplicative $S(ab)=S(b) S(a)$ for all $a, b\in A$.
   \item[(3)] $S$ is anticomultiplicative $\Delta S(a)=S(a_{(2)})\otimes S(a_{(1)})$ for all $a\in A$.
 \end{enumerate}
 Hence a Hopf coquasigroup is Hopf algebra iff it is coassociative.

\subsection{Multiplier algebras}

 Let $A$ be an (associative) algebra. We do not assume that $A$ has a unit, but we do require that
 the product, seen as a bilinear form, is non-degenerated. This means that, whenever $a\in A$ and $ab=0$ for all $b\in A$
 or $ba=0$ for all $b\in A$, we must have that $a=0$.
 Then we can consider the multiplier algebra $M(A)$ of $A$.
 Recall from \cite{V94,V98} that $M(A)$ is characterized as the largest algebra with identity containing $A$ as an essential two-sided ideal.
 In particularly, we still have that, whenever $a\in M(A)$ and $ab=0$ for all $b\in A$ or $ba=0$ for all $b\in A$, again $a=0$.
 Furthermore, we consider the tensor algebra $A\otimes A$. It is still non-degenerated and we have its multiplier algebra $M(A\otimes A)$.
 There are natural imbeddings
 $$A\otimes A \subseteq M(A)\otimes M(A) \subseteq M(A\otimes A).$$
 In generally, when $A$ has no identity, these two inclusions are stict.
 If $A$ already has an identity, the product is obviously non-degenerate
 and $M(A)=A$ and $M(A\otimes A) = A\otimes A$. More details about the concept of the multiplier algebra of an algebra, we refer to \cite{V94}.

 Let $A$ and $B$ be non-degenerate algebras, if homomorphism $f: A\longrightarrow M(B)$ is non-degenerated
 (i.e., $f(A)B=B$ and $Bf(A)=B$),
 then has a unique extension to a homomorphism $M(A)\longrightarrow M(B)$, we also denote it $f$.

\section{Integrals on an infinite dimensional Hopf quasigroup}
\def\theequation{\thesection.\arabic{equation}}
\setcounter{equation}{0}

 Let $H$ be a finite dimensional Hopf quasigroup and $H^*= Hom(H, k)$ be the dual space with natural Hopf coquasigruoup structure given by
 \begin{eqnarray*}
 && \langle ab, h\rangle = \langle a, h_{(1)}\rangle\langle b, h_{(2)}\rangle, \qquad \langle \Delta(a), h\otimes g\rangle = \langle a, hg\rangle, \\
 && \langle 1, h\rangle = \varepsilon(h), \qquad \langle a, 1\rangle = \varepsilon(a), \qquad \langle S(a), h\rangle = \langle a, S(h)\rangle.
 \end{eqnarray*}
 Then there is a natural question: For an infinite dimensional Hopf quasigroup $H$, how about its dual?

 Recall from \cite{K}, a left (resp. right) integral on $H$ is a nonzero element $\varphi \in H^*$ (resp. $\psi \in H^*$ ) such that
 \begin{eqnarray*}
 && (id\otimes \varphi)\Delta(h) = \varphi(h)1_{H}  \quad \big(\mbox{resp.} (\psi\otimes id)\Delta(h) = \psi(h)1_{H} \big), \qquad \forall h\in H.
 \end{eqnarray*}
 And from Lemma 3.3 in \cite{K}, we have that $\varphi\circ S$ is a right integral on $H$.
 \\

 \textbf{Lemma \thesection.1} [\cite{K}, Lemma 3.4, 3.8]
 Let $\varphi$ (resp. $\psi$) be a left (resp. right) integral on $H$, then for $h, g\in H$
 \begin{eqnarray}
 && h_{(1)} \varphi\big( h_{(2)}S(g) \big) = \varphi\big( hS(g_{(1)}) \big) g_{(2)}, \qquad h_{(1)}\varphi (gh_{(2)}) = S (g_{(1)}) \varphi (g_{(2)} h). \label{3.1} \\
 && \psi \big( S(g)h_{(1)}\big) h_{(2)} = \psi \big( S(g_{(2)})h \big) g_{(1)}, \qquad \psi(g_{(1)}h) g_{(2)} = \psi(g h_{(1)}) S(h_{(2)}). \label{3.2}
 \end{eqnarray}

 \emph{Proof} From the proof of Lemma 3.4 and 3.8 in \cite{K}, we can easily check the above equations also hold in infinite dimensional case.
 $\hfill \Box$
 \\

 In the following, we will constuct the 'integral dual' of a class of infinite dimensional Hopf quasigroup.
 Let $H$ be an infinite dimensional Hopf quasigroup with a bijective antipode and a left faithful integral, i.e., $\varphi(gh)=0, \forall h\in H \Rightarrow g=0$ and $\varphi(gh)=0, \forall g\in H \Rightarrow h=0$.

 First, we show that for the infinite dimensional Hopf quasigroup, the left faithful integral  is unique up to scalar.

 \textbf{Theorem \thesection.2}
 Let $\varphi'$ be another left faithful integral on $H$, then $\varphi'=\lambda\varphi$ for some scalar $\lambda\in k$, i.e.,
 the left faithful integral on $H$ is unique up to scalar.

 \emph{Proof}
 From Lemma \thesection.1, we have $h_{(1)}\varphi (gh_{(2)}) = S(g_{(1)}) \varphi (g_{(2)} h)$ for all $h, g\in H$.
 Apply $\varphi'$ to both expressions in this equation. Because $\varphi' \circ S$ is a right integral, the right hand side will give
 \begin{eqnarray*}
 && \varphi'(S g_{(1)}) \varphi (g_{(2)} h) = \varphi (\varphi' S (g_{(1)})g_{(2)} h) = \varphi (\varphi' S (g) 1_{H}\cdot h) = \varphi' S (g) \varphi (h).
 \end{eqnarray*}
 For the left hand side,
 \begin{eqnarray*}
 && \varphi'(h_{(1)}\varphi (gh_{(2)})) = \varphi'(h_{(1)})\varphi (gh_{(2)}) = \varphi (g\varphi'(h_{(1)})h_{(2)}) = \varphi (g \delta_{h}),
 \end{eqnarray*}
 where $\delta_{h} = \varphi'(h_{(1)})h_{(2)}$.
 Therefore, $\varphi' S (g) \varphi (h) = \varphi (g \delta_{h})$ for all $h, g\in H$.

 We claim that there is an element $\delta\in H$ such that $\delta_{h} = \varphi(h)\delta$ for all $h\in H$.
 Indeed, for any $h'\in H$
 \begin{eqnarray*}
 \varphi (g \varphi (h')\delta_{h})
 &=& \varphi (h') \varphi (g \delta_{h}) = \varphi (h') \varphi'S(g) \varphi (h)\\
 &=& \varphi (h) \varphi'S(g) \varphi (h') = \varphi (h) \varphi (g \delta_{h'})\\
 &=& \varphi (g \varphi (h)\delta_{h'}),
 \end{eqnarray*}
 then $\varphi (h')\delta_{h} = \varphi (h)\delta_{h'}$ for all $h, h'\in H$, since $\varphi$ is faithful.
 Choose an $h'\in H$ such that $\varphi (h') = 1$ and denote $\delta = \delta_{h'}$, then $\delta_{h} = \varphi(h)\delta$.

 If we apply $\varepsilon$, we get
 \begin{eqnarray*}
 \varphi(h)\varepsilon(\delta) &=& \varepsilon(\delta_{h}) =  \varepsilon(\varphi'(h_{(1)})h_{(2)}) \\
 &=&  \varphi'(h_{(1)})\varepsilon(h_{(2)}) =  \varphi'(h_{(1)}\varepsilon(h_{(2)}) ) \\
 &=&  \varphi'(h)
 \end{eqnarray*}
 for all $h\in H$ and with $\lambda = \varepsilon(\delta)$, we find the desired result.
 $\hfill \Box$
 \\

 \textbf{Remark}
  Similarly, the right faithful integral on $H$ is unique up to scalar. However, it is a pity that the non-zero faithful integrals do not always exist in infinite dimensional case,
  even for the special infinite dimensional Hopf algebra case.
 \\

 \textbf{Proposition \thesection.3}
 There is a unique group-like element $\delta\in H$ such that for all $h\in H$
 \begin{enumerate}
 \item[(1)] $(\varphi\otimes id)\Delta(h) = \varphi(h)\delta$.
 \item[(2)] $\varphi S(a) = \varphi(a\delta)$.
  \end{enumerate}
 Furthermore, if the antipode $S$ is bijective, then $(id\otimes \psi)\Delta(h) = \psi(h)\delta^{-1}$.

 \emph{Proof}
 From the proof of Proposition \thesection.2, $\varphi(h)\delta = \delta_{h} = \varphi'(h_{(1)})h_{(2)}$ and $\varphi' S (g) \varphi (h) = \varphi (g \delta_{h})$,
 we take $\varphi' = \varphi$ and get a element $\delta\in H$ such that $(\varphi\otimes id)\Delta(h) = \varphi(h)\delta$ and $\varphi S (h) = \varphi (h\delta)$.
 This gives the first part of (1) and (2).

 If we apply $\varepsilon$ and $\Delta$ on the first equation, we find $\varepsilon(\delta)=1$ and $\Delta(\delta) = \delta \otimes \delta$.
 by Proposition 4.2 (1) in \cite{KM} $S(\delta)\delta = 1 =\delta S(\delta)$, then $S(\delta) = \delta^{-1}$. Hence $\delta$ is a group-like element.

 Because $S$ flips the coproduct and if we let $\psi=\varphi\circ S$, we get
 \begin{eqnarray*}
 (id\otimes \psi)\Delta(h)
 &=& S^{-1}(S\otimes \psi)\Delta(h) = S^{-1}(S\otimes \varphi\circ S)\Delta(h) \\
 &=& S^{-1}(id\otimes \varphi) (S\otimes S)\Delta(h) = S^{-1}(id\otimes \varphi) \Delta^{cop}(S(h)) \\
 &=& S^{-1}(\varphi\otimes id) \Delta(S(h)) \stackrel{(1)}{=} S^{-1}(\varphi(S(h))\delta) \\
 &=& \psi(h)\delta^{-1}.
 \end{eqnarray*}
 This completes the proof.
 $\hfill \Box$
 \\

  \textbf{Remark} (1) The square $S^{2}$ leaves the coproduct invariant, it follows that the composition $\varphi\circ S^{2}$ of the left faithful integral $\varphi$ with $S^{2}$ will again a left faithful integral.
  By the uniqueness of left faithful integrals, there is a number $\tau \in k$ such that $\varphi\circ S^{2} = \tau\varphi$.

  (2) If we apply (2) in Proposition \thesection.3 twice, we get
 \begin{eqnarray*}
 && \varphi\big( S^{2}(a) \big) = \varphi\big( S(a)\delta \big) = \varphi\big( S(\delta^{-1} a) \big)  = \varphi\big((\delta^{-1} a) \delta \big).
 \end{eqnarray*}
 So $\varphi\big((\delta^{-1} a) \delta \big) = \tau\varphi(a)$.

\section{Integral dual}
\def\theequation{\thesection.\arabic{equation}}
\setcounter{equation}{0}

 In this section, we will construct the dual of an infinite dimensional Hopf quasigroup.
 This construction bases on the faithful integrals introduced in the last section.
 Here, we also start with defining the following subspace of the dual space $H^*$.

 \textbf{Definition \thesection.1}
 Let $\varphi$ be a left faithful integral on a Hopf quasigroup $H$.
 We define $\widehat{H}$ as the space of linear functionals on $H$ of the form $\varphi(\cdot h)$ where $h\in H$, i.e.,
 \begin{eqnarray*}
 \widehat{H} = \{\varphi(\cdot h) \, | \, h\in H\}.
 \end{eqnarray*}

 \textbf{Lemma \thesection.2}
 Let $H$ be a Hopf quasigroup and $\varphi$ (resp. $\psi$) be a left (resp. right) integral on $H$.
 If $a\in H$, then there is a $b\in H$ such that $\varphi(ax)=\psi(xb)$ for all $x\in H$.
 Similarly, given $q\in H$, we have $p\in H$ so that $\varphi(xp)=\psi(q x)$ for all $x\in H$.

 \emph{Proof}
 By the equations (\ref{3.1}) in Lemma \thesection.1 and (\ref{3.2}), we have for any $h, x\in H$,
  \begin{eqnarray*}
 (\psi\otimes \varphi)\big(xq_{(1)} \otimes pS(q_{(2)})\big)
 &=& \psi(xq_{(1)}) \varphi(pS(q_{(2)})) = \varphi\big(p \underline{\psi(xq_{(1)})S(q_{(2)})} \big) \\
 &\stackrel{(\ref{3.2})}{=}& \varphi\big(p \psi(x_{(1)}q) x_{(2)} \big) =\psi(x_{(1)}q) \varphi\big(p x_{(2)} \big) \\
 &=& \psi\big(\underline{x_{(1)}\varphi(p x_{(2)})} q \big)  \stackrel{(\ref{3.1})}{=} \psi\big(S(p_{(1)})\varphi(p_{(2)} x) q \big) \\
 &=&\varphi(p_{(2)} x) \psi\big(S(p_{(1)}) q \big) = \varphi\Big((\psi\big(S(p_{(1)}) q \big)p_{(2)}) x \Big) \\
 &=& \varphi\Big((\psi\circ S \otimes id)\big( (S^{-1}(q) \otimes 1)\Delta(p) \big) x \Big).
 \end{eqnarray*}
 On the other hand, we also have
  \begin{eqnarray*}
 (\psi\otimes \varphi)\big(xq_{(1)} \otimes pS(q_{(2)})\big)
 &=& \psi(xq_{(1)}) \varphi(pS(q_{(2)})) = \psi\Big(x (q_{(1)}\varphi\big(pS(q_{(2)})\big)) \Big) \\
 &=& \psi\Big(x (id\otimes \varphi\circ S)\big(\Delta(q) (1\otimes S^{-1}(p))\big) \Big).
 \end{eqnarray*}
 By Theorem 4.5 in \cite{WW19}, Galois maps $T_{1}: a\otimes b\mapsto \Delta(a)(1 \otimes b) $ and $T_{2}:  a\otimes b\mapsto (a\otimes 1)\Delta(b)$ are bijective,
 then any element in $H$ has the form $(\psi\circ S \otimes id)\big( (S^{-1}(q) \otimes 1)\Delta(p) \big)$.
 Hence the above calculation will give us the formula $\varphi(ax)=\psi(xb)$ for all $x\in H$.

 Similarly by computing $(\psi\otimes \varphi)\big(S(q_{(2)})x \otimes q_{(1)}S(p) \big)$, we get the second assertion.
 $\hfill \Box$
 \\

 \textbf{Remark}
 (1) In the proof of second part, we need Galois maps $T_{3}:  a\otimes b\mapsto \Delta(a)(b \otimes 1)$ and $T_{4}: a\otimes b\mapsto (1\otimes a)\Delta(b)$ are bijective,
 which follows the fact that the antipode $S$ is bijective, and $T_{3}^{-1}:  a\otimes b\mapsto b_{(2)} \otimes S^{-1}(b_{(1)})a$ and $T_{4}: a\otimes b\mapsto bS^{-1}(a_{(2)}) \otimes a_{(1)}$.

 (2) From Lemma \thesection.2, we get that
 \begin{eqnarray*}
 && \widehat{H} = \{\varphi(\cdot h) \, | \, h\in H\} = \{\psi( h\cdot) \, | \, h\in H\}.\\
 &\mbox{and}& \{\varphi( h\cdot) \, | \, h\in H\} = \{\psi(\cdot h) \, | \, h\in H\}.
 \end{eqnarray*}

 In the following, we need the following assumption to construct the dual.

 \textbf{Assumption \thesection.3}
 \begin{eqnarray}
 && \{\varphi(\cdot h) \, | \, h\in H\} = \{\varphi( h\cdot) \, | \, h\in H\}, \\
 && \varphi\big((\cdot h)h'\big), \varphi\big(h'(h\cdot)\big) \in \widehat{H}, \quad \forall  h, h'\in H.
 \end{eqnarray}

 Following this assumption, we have
 \begin{eqnarray*}
 && \widehat{H} = \{\varphi(\cdot h) \, | \, h\in H\} = \{\psi( h\cdot) \, | \, h\in H\} =  \{\varphi( h\cdot) \, | \, h\in H\} = \{\psi(\cdot h) \, | \, h\in H\}
 \end{eqnarray*}
 and $\psi\big((\cdot h)h'\big),  \psi\big(h'(h\cdot)\big)\in \widehat{H}$.
 \\

 We start by making $\widehat{H}$ into an algebra by dualizing the coproduct.

 \textbf{Proposition \thesection.4}
 For $w,w'\in \widehat{H}$, we can define a linear functional $ww'$ on $H$ by the formula
 \begin{eqnarray}
 && (ww')(h) = (w\otimes w')\Delta(h), \quad \forall h\in H. \label{3.3}
 \end{eqnarray}
 Then $ww'\in \widehat{H}$. This product on $\widehat{H}$ ia associative and non-degenerate.

 \emph{Proof}
 Let $w, w'\in \widehat{H}$ and assume that $w' = \varphi(\cdot m)$ with $m\in H$. we have
 \begin{eqnarray*}
 (ww')(h) &=& (w\otimes \varphi(\cdot m))\Delta(h) = (w\otimes \varphi)\big(\Delta(h)(1\otimes m) \big) \\
 &=& w\big( h_{(1)} \varphi(h_{(2)}m) \big) \stackrel{(\ref{3.1})}{=} w \Big(S^{-1}\big( m_{(1)} \varphi(h m_{(2)}) \big)\Big) \\
 &=& \varphi\Big(h \big(wS^{-1}(m_{(1)}) m_{(2)}\big)\Big)
 \end{eqnarray*}
 We see that the product $ww'$ is well-defined as alinear functional on $H$ and it has the form $\varphi(\cdot g)$, where $g=wS^{-1}(m_{(1)}) m_{(2)}$.
 So $ww' \in \widehat{H}$. Therefore, we have defined a product in $\widehat{H}$.

 The associativity of this product in $\widehat{H}$ is an consequence of the coassociativity of $\Delta$ on $H$.

 To prove that the product is non-degenerate, assume that $ww' =0$ for all $w\in \widehat{H}$.
 From the above calculation, for any $h\in H$, $0= (ww')(h) = \varphi\Big(h \big(wS^{-1}(m_{(1)}) m_{(2)}\big)\Big)$,
 then $wS^{-1}(m_{(1)}) m_{(2)} = 0$ because of the faithfulness of $\varphi$.
 This implies $wS^{-1}(m) = 0$ for all $w \in \widehat{H}$, i.e., $\varphi(S^{-1}(m)h) = 0$ for all $h \in H$.
 We conclude that $S^{-1}(m)=0$ then $m=0$, i.e., $w' =0$.
 Similarly, $ww' =0$ for all $w'\in \widehat{H}$ implies $w=0$.
 $\hfill \Box$
 \\

 \textbf{Remark}
 Under the assumption, the elements of $\widehat{H}$ can be expressed in four different forms.
 When we use these different forms in the definition of product in $\widehat{H}$, we get the following useful expressions:
 \begin{enumerate}
 \item[(1)] $w\varphi(\cdot a) = \varphi(\cdot b)$  \mbox{with}  $b = wS^{-1}(a_{(1)}) a_{(2)}$;
 (2) $w\varphi( a\cdot) = \varphi(c \cdot)$  \mbox{with} $c = wS(a_{(1)}) a_{(2)}$.
 \item[(3)] $\psi(\cdot a)w = \psi(\cdot d)$  \mbox{with}  $d = a_{(1)} wS(a_{(2)})$; \quad
 (4)  $\psi( a\cdot)w = \psi(e \cdot)$  \mbox{with}  $e = a_{(1)} wS^{-1}(a_{(2)})$.
 \end{enumerate}

 Moreover, the multiplier algebra $M(\widehat{H})$ of $\widehat{H}$ can be identified with the space $H^*$.
 Indeed, for $f\in H^*$ and $w\in \widehat{H}$, $fw, wf\in \widehat{H}$;
 he counit $\varepsilon$, as a linear functional on $H$, is in fact the unit in the multiplier algebra $M(\widehat{H})$;
 $fw=0$ (resp. $wf=0$) for all $w\in \widehat{H}$ implies $f=0$.
 \\

 Let us now define the comultiplication $\widehat{\Delta}$ on $\widehat{H}$.
 Roughly speaking, the coproduct is dual to the multiplication in $H$ in the sense that
 \begin{eqnarray*}
 \langle\widehat{\Delta}(w), x\otimes y\rangle = \langle w, xy\rangle, \quad \forall x, y\in H.
 \end{eqnarray*}

 \textbf{Definition \thesection.5}
 Let $w_{1}, w_{2}\in \widehat{H}$, then we put
 \begin{eqnarray}
 \langle(w_{1}\otimes 1)\widehat{\Delta}(w_{2}), x\otimes y\rangle &=& \langle w_{1}\otimes w_{2}, x_{(1)}\otimes x_{(2)} y \rangle \label{3.4} \\
 \langle\widehat{\Delta}(w_{1})(1 \otimes w_{2}), x\otimes y\rangle &=& \langle w_{1}\otimes w_{2}, xy_{(1)}\otimes y_{(2)} \rangle \label{3.5}
 \end{eqnarray}
 for all $x, y\in H$.
 \\

 We will first show that the functionals in Definition \thesection.7 are well-defined and again in $\widehat{H}\otimes \widehat{H}$.

 \textbf{Lemma \thesection.6} $(w_{1}\otimes 1)\widehat{\Delta}(w_{2}), \widehat{\Delta}(w_{1})(1 \otimes w_{2}) \in \widehat{H}\otimes \widehat{H}$.
 These above two formulas define $\widehat{\Delta}(w)$ as a multiplier in $M(\widehat{H}\otimes \widehat{H})$ for all $w\in \widehat{H}$.

 \emph{Proof} Let $w_{1} = \psi(a\cdot)$ and $w_{2} = \psi(b\cdot)$, where $a, b\in H$. For any $x, y\in H$, we have
 \begin{eqnarray*}
 \langle(w_{1}\otimes 1)\widehat{\Delta}(w_{2}), x\otimes y\rangle
 &=& \langle w_{1}\otimes w_{2}, x_{(1)}\otimes x_{(2)} y \rangle \\
 &=& \psi(ax_{(1)}) \psi\big(b (x_{(2)} y) \big) = \psi\big(b (\psi(ax_{(1)})x_{(2)} y) \big) \\
 &\stackrel{(\ref{3.2})}{=}& \psi\big(b (\psi(a_{(1)}x)S^{-1}(a_{(2)}) y) \big) = \psi(a_{(1)}x) \psi\big(b (S^{-1}(a_{(2)}) y) \big) \\
 &=& \Big(\psi(a_{(1)} \cdot ) \otimes \psi\big(b (S^{-1}(a_{(2)}) \cdot)\big) \Big) (x\otimes y).
 \end{eqnarray*}
 By the assumption, we obtain that $(w_{1}\otimes 1)\widehat{\Delta}(w_{2})$ is a well-defined element in $\widehat{H}\otimes \widehat{H}$.
 It is similar for the $\widehat{\Delta}(w_{1})(1 \otimes w_{2})$.

 Using the fact that the product in $\widehat{H}$ is dual to the coproduct in $H$ and $\Delta$ in $H$ is coassociative,
 it easily followings that $\big((w_{1}\otimes 1)\widehat{\Delta}(w_{2})\big)(1\otimes w_{3}) = (w_{1}\otimes 1)\big(\widehat{\Delta}(w_{2})(1\otimes w_{3})  \big) $.
 Therefore, $\widehat{\Delta}(w)$ is defined as a two-side  multiplier in $M(\widehat{H}\otimes \widehat{H})$.
 $\hfill \Box$
 \\

 \textbf{Proposition \thesection.7}
 $\widehat{\Delta}: \widehat{H} \longrightarrow M(\widehat{H}\otimes \widehat{H})$ is an algebra homomorphism,
 and also  $(1 \otimes w_{1})\widehat{\Delta}(w_{2}), \widehat{\Delta}(w_{1})(w_{2} \otimes 1) \in \widehat{H}\otimes \widehat{H}$.

 \emph{Proof}
 It is straightforward that $\widehat{\Delta}$ is an algebra homomorphism, since for all $x, y\in H$
 \begin{eqnarray*}
 \langle\widehat{\Delta}(w_{1}w_{2})(1\otimes w_{3}), x\otimes y\rangle &=& \langle w_{1}w_{2}\otimes  w_{3}, xy_{(1)} \otimes y_{(2)} \rangle  \\
 &=& \langle w_{1}, x_{(1)}y_{(1)(1)} \rangle\langle w_{2}, x_{(2)}y_{(1)(2)} \rangle \langle w_{3}, y_{(2)} \rangle, \\
 \langle\widehat{\Delta}(w_{1})\widehat{\Delta}(w_{2}) (1\otimes w_{3}), x\otimes y\rangle
 &=&  \langle\widehat{\Delta}(w_{1}) (f\otimes g) , x\otimes y \rangle \qquad \mbox{($\widehat{\Delta}(w_{2}) (1\otimes w_{3}) := f\otimes g$ )}\\
 &=& \langle\widehat{\Delta}(w_{1})(1\otimes g), x_{(1)}\otimes y \rangle \langle f, x_{(2)} \rangle \\
 &=& \langle w_{1} \otimes g, x_{(1)}y_{(1)} \otimes y_{(2)} \rangle \langle f, x_{(2)} \rangle \\
 &=& \langle w_{1}, x_{(1)}y_{(1)}\rangle \langle f \otimes g, x_{(2)} \otimes y_{(2)} \rangle \\
 &=& \langle w_{1}, x_{(1)}y_{(1)}\rangle \langle \widehat{\Delta}(w_{2}) (1\otimes w_{3}), x_{(2)} \otimes y_{(2)} \rangle \\
 &=& \langle w_{1}, x_{(1)}y_{(1)}\rangle \langle w_{2}\otimes w_{3}, x_{(2)}y_{(2)(1)} \otimes y_{(2)(2)} \rangle \\
 &=& \langle w_{1}, x_{(1)}y_{(1)}\rangle \langle w_{2}, x_{(2)}y_{(2)(1)}\rangle \langle w_{3}, y_{(2)(2)} \rangle.
 \end{eqnarray*}
 By the coassociativity of $\Delta$ of $H$, 
 we get $\widehat{\Delta}(w_{1}w_{2})(1\otimes w_{3}) = \widehat{\Delta}(w_{1})\widehat{\Delta}(w_{2}) (1\otimes w_{3})$ for all $w_{3} \in \widehat{H}$.
 This implies $\widehat{\Delta}(w_{1}w_{2}) = \widehat{\Delta}(w_{1})\widehat{\Delta}(w_{2})$.
 
 With the bijective antipode, the proof of the second assertion is similar to the proof of Lemma \thesection.8.
 $\hfill \Box$
 \\

 Let $w\in \widehat{H}$ and assume $w=\varphi(\cdot a)$ with $a\in H$ then. Define $\widehat{\varepsilon}(w)=\varphi(a)=w(1_{H})$.
 Then $\widehat{\varepsilon}$ is a counit on $(\widehat{H}, \widehat{\Delta})$ as follows.

 \textbf{Proposition \thesection.8}
 $\widehat{\varepsilon}: \widehat{H}\longrightarrow k$ is an algebra homomorphism satisfying
 \begin{eqnarray}
 (id\otimes \widehat{\varepsilon})\big((w_{1}\otimes 1)\widehat{\Delta}(w_{2}) \big) =  w_{1}w_{2} \\
 (\widehat{\varepsilon}\otimes id)\big( \widehat{\Delta}(w_{1})(1 \otimes w_{2})\big) =  w_{1}w_{2}
 \end{eqnarray}
 for all $w_{1}, w_{2}\in \widehat{H}$.

 \emph{Proof}
 Firstly, let $w_{1} = \varphi(a\cdot)$ and $w_{2} = \varphi(b\cdot)$, then
 $w_{1}w_{2} = \varphi(c\cdot)$ with $c=\varphi\big( aS(b_{(1)})\big) b_{(2)}$.
 Therefore, if $\psi = \varphi \circ S$ we have
 \begin{eqnarray*}
 \widehat{\varepsilon}(w_{1}w_{2}) &=& \varphi(c) = \varphi\big( aS(b_{(1)})\big) \varphi(b_{(2)}) \\
 &=& \varphi\big( aS(b_{(1)}\varphi(b_{(2)}))\big) = \varphi(a) \varphi(b) \\
 &=& \widehat{\varepsilon}(w_{1})\widehat{\varepsilon}(w_{2}).
 \end{eqnarray*}

 Secondly, let $w_{1} = \psi(a\cdot)$ and $w_{2} = \psi(b\cdot)$, then we have
 \begin{eqnarray*}
 (w_{1}\otimes 1)\widehat{\Delta}(w_{2}) &=& \psi(a_{(1)} \cdot ) \otimes \psi\big(b (S^{-1}(a_{(2)}) \cdot) \big), \\
 \psi(a\cdot)\psi(b\cdot) &=& \psi(a_{(1)} \cdot ) \psi\big(b S^{-1}(a_{(2)}) \big).
 \end{eqnarray*}
 Hence,
 \begin{eqnarray*}
 (id\otimes \widehat{\varepsilon})\big((w_{1}\otimes 1)\widehat{\Delta}(w_{2}) \big)
 &=& \psi(a_{(1)} \cdot ) \psi\big(b (S^{-1}(a_{(2)}) 1) \big) \\
 &=& \psi(a_{(1)} \cdot ) \psi\big(b S^{-1}(a_{(2)}) \big)  \\
 &=& w_{1} w_{2}.
 \end{eqnarray*}

 Finally, the second formula is proven in a similar way, in this case letting $w_{1} = \varphi(\cdot a)$ and $w_{2} = \varphi(\cdot b)$.
 $\hfill \Box$
 \\

 Let $\widehat{S}: \widehat{H}\longrightarrow \widehat{H}$ be the dual to the antipode of $H$, i.e., $\widehat{S}(w)=w\circ S$.
 Then it is easy to see that $\widehat{S}(w) \in \widehat{H}$, and we have the following property.

  \textbf{Proposition \thesection.9}
  $\widehat{S}$ is antimultiplicative and coantimultiplicative such that
  \begin{eqnarray*}
  w'\otimes w
  &=& (m\otimes id)(id \otimes \widehat{S}\otimes id)\Big((id\otimes\widehat{\Delta})\big( (w'\otimes 1)\widehat{\Delta}(w)\big) \Big) \\
  &=& (m\otimes id)(\widehat{S}\otimes id\otimes id)\Big((id\otimes \widehat{\Delta})\big( \widehat{\Delta}(w) (\widehat{S}^{-1}(w')\otimes 1)\big)\Big) \\
  &=& (id\otimes m)(id \otimes \widehat{S}\otimes id)\Big((\widehat{\Delta}\otimes id)\big( \widehat{\Delta}(w') (1\otimes  w)\big) \Big) \\
  &=& (id\otimes m)(id \otimes id\otimes\widehat{S})\Big((\widehat{\Delta}\otimes id)\big( (1\otimes \widehat{S}^{-1}(w))\widehat{\Delta}(w')\big)\Big).
 \end{eqnarray*}

 \emph{Proof}
 For $w_{1}, w_{2}\in \widehat{H}$ and any $x\in H$,
 \begin{eqnarray*}
 \langle\widehat{S}(w_{1}w_{2}), x\rangle &=& \langle w_{1}w_{2}, S(x)\rangle = \langle w_{1}, S(x_{2})\rangle\langle w_{2}, S(x_{1})\rangle \\
 &=& \langle \widehat{S}(w_{1}), x_{2}\rangle\langle \widehat{S}(w_{2}), x_{1}\rangle = \langle \widehat{S}(w_{2})\widehat{S}(w_{1}), x\rangle
 \end{eqnarray*}
 This implies $\widehat{S}$ is antimultiplicative.
 \begin{eqnarray*}
 \langle\widehat{\Delta}\widehat{S}(w_{1})(1\otimes S(w_{2})), x \otimes y\rangle
 &\stackrel{(\ref{3.5})}{=}& \langle\widehat{S}(w_{1})\otimes \widehat{S}(w_{2}), xy_{(1)} \otimes y_{(2)}\rangle \\
 &=& \langle w_{1}, S(xy_{(1)}) \rangle \langle w_{2}, S(y_{(2)})\rangle \\
 &=& \langle w_{1}, S(y)_{(2)}S(x) \rangle \langle w_{2}, S(y)_{(1)}\rangle \\
 &=& \langle w_{2}\otimes w_{1}, S(y)_{(1)}\otimes S(y)_{(2)}S(x) \rangle \\
 &\stackrel{(\ref{3.4})}{=}& \langle (w_{2}\otimes 1)\widehat{\Delta}(w_{1}), S(y)\otimes S(x) \rangle \\
 &=& \langle(\widehat{S}\otimes\widehat{S}) \widehat{\Delta}^{cop}(w_{1})(1\otimes S(w_{2})), x \otimes y\rangle,
 \end{eqnarray*}
 We conclude $\widehat{S}$ is coantimultiplicative.

 Finally, we show $w'\otimes w  = (m\otimes id)(id \otimes \widehat{S}\otimes id)\Big((id\otimes\widehat{\Delta})\big( (w'\otimes 1)\widehat{\Delta}(w)\big) \Big)$,
  the other three formulas is similar.
 \begin{eqnarray*}
 && \langle (m\otimes id)(id \otimes \widehat{S}\otimes id)\Big((id\otimes\widehat{\Delta})\big( (w'\otimes 1)\widehat{\Delta}(w)\big) \Big), x \otimes y\rangle \\
 &\stackrel{(\ref{3.3})}{=}& \langle (id\otimes\widehat{\Delta})\big( (w'\otimes 1)\widehat{\Delta}(w)\big) , x_{(1)} \otimes S(x_{(2)}) \otimes y\rangle \\
 &=& \langle \big( (w'\otimes 1)\widehat{\Delta}(w)\big), x_{(1)} \otimes S(x_{(2)})y \rangle \\
 &\stackrel{(\ref{3.4})}{=}& \langle  w'\otimes w, x_{(1)(1)} \otimes x_{(1)(2)}(S(x_{(2)})y) \rangle \\
 &=& \langle  w'\otimes w, x \otimes y \rangle.
 \end{eqnarray*}
 This completes the proof.
 $\hfill \Box$
 \\

 The equation in the Proposition \thesection.11 can be expressed by generalized Sweedler notation as follows.
  \begin{eqnarray*}
  w'\otimes w
  &=& w' \widehat{S}(w_{(1)}) w_{(2)(1)} \otimes w_{(2)(2)} = w'w_{(1)} \widehat{S}(w_{(2)(1)}) \otimes w_{(2)(2)}\\
  &=& w'_{(1)(1)} \otimes \widehat{S}(w'_{(1)(2)}) w'_{(2)} w = w'_{(1)(1)} \otimes w'_{(1)(2)} \widehat{S}(w'_{(2)}) w.
 \end{eqnarray*}
 As a consequence, the antipode $\widehat{S}$ also satisfies
 \begin{eqnarray*}
 && m(id\otimes \widehat{S})\big( (w_{1}\otimes 1)\widehat{\Delta}(w_{2}) \big) = \widehat{\varepsilon}(w_{2})w_{1}, \\
 && m(\widehat{S}\otimes id)\big( \widehat{\Delta}(w_{1})(1 \otimes w_{2}) \big) = \widehat{\varepsilon}(w_{1})w_{2}.
 \end{eqnarray*}
 In fact, there is another way to prove.
 \begin{eqnarray*}
 \langle m(id\otimes \widehat{S})\big( (w_{1}\otimes 1)\widehat{\Delta}(w_{2}) \big), x \rangle
 &=& \langle (id\otimes \widehat{S})\big( (w_{1}\otimes 1)\widehat{\Delta}(w_{2}) \big), x_{(1)} \otimes x_{(2)}\rangle \\
 &=& \langle \big( (w_{1}\otimes 1)\widehat{\Delta}(w_{2}) \big), x_{(1)} \otimes S(x_{(2)})\rangle \\
 &\stackrel{(\ref{3.4})}{=}& \langle  w_{1}\otimes w_{2}) , x_{(1)(1)} \otimes x_{(1)(2)}S(x_{(2)})\rangle  \\
 &=& \langle  w_{1}\otimes w_{2} , x \otimes 1 \rangle  \\
 &=& \widehat{\varepsilon}(w_{2})w_{1}.
 \end{eqnarray*}

 Let $\psi$ be a right faithful integral on $H$. For $w=\psi(a\cdot)$ we set $\widehat{\varphi}(w) = \varepsilon(a)$. Then we have the following result.

 \textbf{Proposition \thesection.10}
 $\widehat{\varphi}$ defined above is a left faithful integral on $\widehat{H}$.

 \emph{Proof}
 It is clear that $\widehat{\varphi}$ is non-zero. Assume $w_{1} = \psi(a\cdot)$ and $w_{2} = \psi(b\cdot)$ with $a, b\in H$, then
 \begin{eqnarray*}
 (w_{1}\otimes 1)\widehat{\Delta}(w_{2}) &=& \psi(a_{(1)} \cdot ) \otimes \psi\big(b (S^{-1}(a_{(2)}) \cdot) \big).
 \end{eqnarray*}
 Therefore, we have
 \begin{eqnarray*}
 (id\otimes \widehat{\varphi})\Big((w_{1}\otimes 1)\widehat{\Delta}(w_{2})\Big)
 &=& \psi(a_{(1)} \cdot ) \otimes  \widehat{\varphi} \psi\big(b (S^{-1}(a_{(2)}) \cdot) \big) \\
 &=& \psi(a \cdot ) \varepsilon(b)
 = \widehat{\varphi}(w_{2}) w_{1}.
 \end{eqnarray*}

 Next, we show that $\widehat{\varphi}$  is  faithful.
 If $w_{1}, w_{2}\in \widehat{H}$ and assume $w_{1} = \psi(a\cdot)$ with $a\in A$, we have
 $w_{1} w_{2} =\psi \big(a_{(1)} w_{2}S^{-1}(a_{(2)}) \cdot \big)$.
 Therefore, $\widehat{\varphi}(w_{1} w_{2}) = w_{2}S^{-1}(a)$.
 If this is 0 for all $a\in H$, then $w_{2} = 0$,
 wile if this is 0 for all $ w_{2}$ then $a=0$.
 This proves the faithfulness of $\widehat{\varphi}$.
 $\hfill \Box$
 \\

 Now, we introduce an algebraic structure: multipler Hopf coquasigroup, generalizing the ordinary Hopf coquasigroup to a nonunital case.
 Let $A$ be an (associative) algebra, may not has a unit, but the product, seen as a bilinear form, is non-degenerated.
 \\

 \textbf{Definition \thesection.11}
  A \emph{multipler Hopf coquasigroup} is a nondegenerate associative algebra $A$ equipped with algebra homomorphisms $\Delta: A\longrightarrow M(A\otimes A)$(coproduct),
 $\varepsilon: A\longrightarrow k$(counit) and a linear map $S: A \longrightarrow A$ (antipode) such that
 \begin{enumerate}
 \item[(1)] $T_{1}(a\otimes b)=\Delta(a)(1\otimes b)$ and $T_{2}(a\otimes b)=(a\otimes 1)\Delta(b)$ belong to $A\otimes A$ for any $a, b\in A$.
 \item[(2)] The counit satisfies $(\varepsilon\otimes id)T_{1}(a\otimes b) = ab = (id\otimes \varepsilon)T_{2}(a\otimes b)$.
 \item[(3)] $S$ is antimultiplicative and anticomultiplicative such that for any $a, b\in A$
 \begin{eqnarray}
 && S(a_{(1)}) a_{(2)(1)} \otimes a_{(2)(2)} = 1\otimes a = a_{(1)} S(a_{(2)(1)}) \otimes a_{(2)(2)}, \label{4.8} \\
 && a_{(1)(1)} \otimes  S(a_{(1)(2)}) a_{(2)} = a\otimes 1 = a_{(1)(1)} \otimes a_{(1)(2)} S(a_{(2)}). \label{4.9}
 \end{eqnarray}
 \end{enumerate}
 If the antipode $S$ is bijective, then  multipler Hopf coquasigroup $(A, \Delta)$ is called \emph{regular}.
 \\

 \textbf{Remark} 
 (1) In multipler Hopf coquasigroup $(A, \Delta)$, $T_{1}$ and $T_{2}$ are bijective.
 If $(A, \Delta)$ is regular, then $T_{3}$ and $T_{4}$ are also.
 In fact, from (3) in Definition \thesection.11 we can easily get
 \begin{eqnarray*}
 && m(id\otimes {S})\big( (a\otimes 1){\Delta}(b) \big) = {\varepsilon}(b)a, \\
 && m({S}\otimes id)\big( {\Delta}(a)(1 \otimes b) \big) = {\varepsilon}(a)b.
 \end{eqnarray*}

 (2) The equation (\ref{4.8}) and (\ref{4.9}) make sense. Take (\ref{4.8}) for example, (\ref{4.9}) is similar.
 \begin{eqnarray*}
 b\otimes a c
 &=& ba_{(1)} S(a_{(2)(1)}) \otimes a_{(2)(2)} c\\
 &=& (m\otimes id)(id \otimes S\otimes id)\Big((id\otimes \Delta)\big( (b\otimes 1) \Delta(a)\big)(1\otimes 1\otimes c) \Big).
 \end{eqnarray*}
 $ba_{(1)} \otimes a_{(2)} = (b\otimes 1) \Delta(a) \in A\otimes A$, 
 and then $a_{(2)(1)} \otimes a_{(2)(2)} c \in A\otimes A$. 
 Therefore, $b\otimes ac = ba_{(1)} S(a_{(2)(1)}) \otimes a_{(2)(2)} c$ holds for all $b, c\in A$.
 This implies $ a_{(1)} S(a_{(2)(1)}) \otimes a_{(2)(2)} = 1\otimes a$.
 \begin{eqnarray*}
 b\otimes ac 
 &=& S(a_{(1)}) a_{(2)(1)}b \otimes a_{(2)(2)}c  = S(a_{(1)}) a_{(2)(1)}x_{(1)} \otimes a_{(2)(2)}x_{(2)}y \\
 &=& (m\otimes id)(S \otimes id \otimes id)\Big((id\otimes \Delta)\big( \Delta(a)(1\otimes x)\big)(1\otimes 1\otimes y) \Big),
 \end{eqnarray*}
 where $b\otimes c = \Delta(x)(1\otimes y)$.
 $b\otimes ac = S(a_{(1)}) a_{(2)(1)}b \otimes a_{(2)(2)}c$ for all $b, c\in A$ implies $1\otimes a = S(a_{(1)}) a_{(2)(1)} \otimes a_{(2)(2)}$.
 
 (3) The comultiplication  may be not coassociative. Multipler Hopf coquasigroup weakens  the coassociativity of coproduct in multiplier Hopf algebra, 
 while algebraic quantum hypergroup in \cite{LV} weakens the homomorphism of coproduct. This is the main difference. 
 \\

 Following Definition \thesection.11, we get the main result of this section.

 \textbf{Theorem \thesection.12}
 Let $(H, \Delta)$ be an infinite dimensional Hopf quasigroup with a faithful integral $\varphi$ and a bijective antipode $S$.
 Then under Assumption \thesection.3 the integral dual $(\widehat{H}, \widehat{\Delta})$ is a regular multipler Hopf coquasigroup with a faithful integral.
 \\
 
 Because (infinite dimensional) Hopf quasigroup $(H, \Delta)$ has the unit $1_{H}$, then there is a special element $\varphi = \varphi(\cdot 1_{H}) \in \widehat{H}$ 
 such that for $w\in \widehat{H}$
 \begin{eqnarray*}
 (w\varphi)(h) = (w\otimes\varphi)\Delta(h) = \varphi(h)w(1) = \widehat{\varepsilon}(w)\varphi(h).
 \end{eqnarray*}
 This implies $w\varphi = \widehat{\varepsilon}(w)\varphi$. 
 We call $\varphi$ a cointegral in $(\widehat{H}, \widehat{\Delta})$. 

 Analogous to multiplier Hopf algebra case in \cite{V98}, we say that a regular multiplier Hopf coquasigroup with a faithful integral $(A, \Delta)$ is of discrete type, if there is a non-zero element $\xi \in A$ so that $a\xi=\varepsilon(a)\xi$ for all $a\in A$.
 
 Then we have the integral dual $(\widehat{H}, \widehat{\Delta})$ of infinite dimensional Hopf quasigroup $(H, \Delta)$
 is a multipler Hopf coquasigroup of discrete type.

\section{Multiplier Hopf coquasigroup: motivating example}
\def\theequation{\thesection.\arabic{equation}}
\setcounter{equation}{0}

 In last section, we introduce the notion of multiplier Hopf coquasigroup, extending Hopf coquasigroup to a nonunital case, 
 and provided a interesting construction: the integral dual of infinite dimensional Hopf quasigroups with integrals.

 In the following, we firstly introduce the motivating example, where Assumption 4.3 naturally holds.
 And then we make some direct comments on multiplier Hopf coquasigroups.
 \\
 
 \textbf{Example \thesection.1}
 Let $G$ be a infinite (IP) quasigroup with identity element $e$, by definition $u^{-1}(uv) = v =(vu)u^{-1}$ for all $u, v\in G$.
 We have that the quasigroup algebra $kG$ is a Hopf quasigroup with the structure shown on the base element $\{u | u\in G\}$
 \begin{eqnarray*}
 \Delta(u) = u\otimes u, \quad \varepsilon(u) = 1, \quad S(u)=u^{-1}.
 \end{eqnarray*}
 The function $\delta_{u}, u\in G$ on $kG$ is given by $\delta_{u}(v) = \delta_{u, v}$, where $\delta_{u, v}$ is the kronecker delta.
 Then $\delta_{e}$ is the left and right integral on $kG$.

 The integral dual $k(G) = \widehat{kG} = \{\delta_{e}(\cdot u) | u\in kG\} = \{\delta_{u^{-1}} | u\in kG\} = \{\delta_{e}(u\cdot) | u\in kG\}$.
 $\delta_{e}\big((\cdot u)v\big) = \delta_{v^{-1} u^{-1}} \in k(G)$ and $\delta_{e}\big( u( v\cdot)\big) = \delta_{v^{-1} u^{-1}} \in k(G)$. 
 Assumption 4.3 naturally holds. Then $(k(G), \widehat{\Delta}, \widehat{\varepsilon}, \widehat{S})$ is a multipler Hopf coquasigroup with the structure as follows.

 As an algebra, $k(G)$ is a nondegenerate algebra with the product
 \begin{eqnarray*}
 \delta_{u}\delta_{v} = \delta_{u, v}\delta_{v},
 \end{eqnarray*}
 and $1 = \sum_{u\in G} \delta_{u}$ is the unit in $M(k(G))$.
 The coproduct, counit and antipode are given by
 \begin{eqnarray*}
 \widehat{\Delta}(\delta_{u}) = \sum_{v\in G}\delta_{v}\otimes \delta_{v^{-1} u}, \quad \widehat{\varepsilon}(\delta_{u}) = \delta_{u, e},\quad \widehat{S}(\delta_{u}) = \delta_{u^{-1}}.
 \end{eqnarray*}

 By the definition of $\widehat{\varphi}$, we get the left integral on $k(G)$ is the function that maps every $\delta_{u}$ to $1$.
 \\ 
 
 As in the theory of multiplier Hopf algebra in \cite{V94}, we also can define a multiplier Hopf $*$-coquasigroup $(A, \Delta)$ over $\mathds{C}$, in which
 $(A, \Delta)$ is a regular multiplier Hopf coquasigroup with the coproduct, counit and antipode compatible with the involution $*$. i.e.,
 \begin{enumerate}
 \item[(1)] The comultiplication $\Delta$ is also a $*$-homomorphism (i.e., $\Delta(a^*) = \Delta(a)^*$);
 \item[(2)] $\varepsilon(a^*) = \overline{\varepsilon(a)}$, where $\overline{(\cdot)}$ means the conjugation of complex numbers;
 \item[(3)] $S(S(a)^*)^* = a$.
 \end{enumerate} 
 
 \textbf{Example \thesection.2}
 In Example \thesection.1 if $k = \mathds{C}$, then $\mathds{C}(G)$ is a multiplier Hopf $*$-coquasigroup.
 \\

 \textbf{Proposition \thesection.3} 
 Let  $(A, \Delta)$ be a multipler Hopf ($*$-)coquasigroup.
 Then $(A, \Delta)$ is a multiplier Hopf ($*$-)algebra introduced, if and only if the comultiplication $\Delta$ is coassociative.
 \\

 \textbf{Proposition \thesection.4}
 If multipler Hopf coquasigroup $(A, \Delta)$ has the unit $1$, then $(A, \Delta)$ is the usual Hopf coquasigroup.
 \\

 Following these two results, multipler Hopf coquasigroup can be considered as the generalization of multiplier Hopf algebra and Hopf coquasigroup.
 Naturally, we can define flexible, alternative and Moufang multipler Hopf coquasigroup.
  
 A multipler Hopf coquasigroup $(A, \Delta)$ is called \emph{flexible} if
 \begin{eqnarray}
 a_{(1)}a_{(2)(2)}\otimes a_{(2)(1)} = a_{(1)(1)}a_{(2)} \otimes a_{(1)(2)}, \quad \forall a\in A,
 \end{eqnarray}
 and \emph{alternative} if also
 \begin{eqnarray}
 && a_{(1)}a_{(2)(1)}\otimes a_{(2)(2)} = a_{(1)(1)}a_{(1)(2)} \otimes a_{(2)}, \\
 && a_{(1)}\otimes a_{(2)(1)}a_{(2)(2)} = a_{(1)(1)} \otimes a_{(1)(2)}a_{(2)}, \quad \forall a\in A.
 \end{eqnarray}
 $A$ is called \emph{Moufang} if
 \begin{eqnarray}
 a_{(1)}a_{(2)(2)(1)}\otimes a_{(2)(1)} \otimes a_{(2)(2)(2)} = a_{(1)(1)(1)}a_{(1)(2)} \otimes a_{(1)(1)(2)}\otimes a_{(2)}, \quad \forall a\in A,
 \end{eqnarray}

 \textbf{Remark}
 (1) By 'cover technique' introduced in \cite{V08}, these four equations make sense.
 
 (2) From the dual, we can get that the integral dual $(\widehat{H}, \widehat{\Delta})$ of infinite dimensional flexible (resp. alternative, Moufang) Hopf quasigroup $(H, \Delta)$
 is a flexible (resp. alternative, Moufang) multipler Hopf coquasigroup.

\section*{Acknowledgements}

 The work was partially supported by the China Postdoctoral Science Foundation (No. 2019M651764)
 and National Natural Science Foundation of China (No. 11601231).


\end {document}